\documentclass[12pt]{article}
\usepackage{latexsym, amsmath, amssymb, amsfonts, amscd}
\usepackage{color}
\usepackage{graphics}
\usepackage[dvips]{epsfig}

\newtheorem{theorem}{Theorem}
\newtheorem{proposition}{Proposition}
\newtheorem{lemma}{Lemma}

\newtheorem{conjecture}{Conjecture}

\newcommand{\adL}{\mbox{\rm ad}_{\Lambda}}

\newcommand{\ad}[1]{\mbox{\rm ad}_{#1}}

\def\lcf{\lbrack\! \lbrack}
\def\rcf{\rbrack\! \rbrack}

\def\dbar{\overline\partial}

\def\oom{\overline\omega}

\newcommand{\ZZ}{\mathbb{Z}}

\def\dbar{\overline\partial}

\def\tri{\triangle}

\def\oomega{\overline\omega}
\def\orho{\overline\rho}
\def\Oomega{\overline\Omega}

\newcommand{\lieg}{\mathfrak{g}}
\newcommand{\liet}{\mathfrak{t}}
\newcommand{\liec}{\mathfrak{c}}

\newcommand{\lbra}[2]{\lcf #1, #2 \rcf}

\newcommand{\Lie}[1]{\operatorname{\textsl{#1}}}

\newcommand{\Sch}{\Lie{Sch}}
\newcommand{\Ger}{\Lie{Ger}}

\newcommand{\image}{\operatorname{\mbox{\rm image} }}

\newcommand{\bproof}{\noindent{\it Proof: }}
\newcommand{\eproof}{\hfill \qed \vspace{0.2in}}
\def\qed{\rule{2.3mm}{2.3mm}}

\begin{document}
\title{\bf
Frobenius Structures \\ and Generalized Deformation \\ of Kodaira Manifolds}
\author{
Yat Sun  Poon\thanks{ Address:
    Department of Mathematics, University of California at Riverside,
    Riverside, CA 92521, U.S.A.. E-mail: ypoon@ucr.edu.} }
\date{September 20, 2020}
\maketitle
\begin{abstract}  It is known that generalized deformation in the sense of Hitchin-Gaultieri is a geometric realization of the degree-$2$ component of 
Kontsevich-Barannikov's homological approach to extended deformation. Through extended deformation, one associates a Frobenius structure
to the extended moduli space. In this notes, we prove that on primary Kodaira manifolds the restriction of the Frobenius structure on the degree-$2$ component
of the extended moduli space is trivial. It generalizes the author's past observation on Kodaira surface. 
\end{abstract}

\section{Introduction}

Inspired by the development of mirror symmetry, especially its homological aspects, 
Kontsevich and Barannikov developed an \it extended \rm deformation 
theory to study mirror symmetry and its higher dimension analogue
 \cite{Barannikov} \cite{B-Kontsevich}.
In \cite{Barannikov}, its starting point is with the complex manifold side of a mirror pair.
 The classical deformation theory of complex structures focuses on the first Dolbeault cohomology
  $H^1(M, \Theta)$  where $M$ is a complex manifold and $\Theta$ 
  is its holomorphic tangent bundle.
  Whereas in an extended deformation theory, Barannikov and Kontsevich utilize 
  the full cohomology with holomorphic polyvector
  fields
  \[
  H^\bullet(M, \wedge^\bullet\Theta)=\sum_{p,q}H^q(M, \wedge^p\Theta)
  \]
  to capture the variation of cohomological structures, 
  e.g. Hodge structure. In this development, one of their discoveries is 
  a Frobenius structure on the extended "moduli" space as
  a $\ZZ$-graded supermanifold
  \cite{B-Kontsevich} \cite{Manin} \cite{Merkulov}. 
  This author had computed the Frobenius structure associated to the primary
  Kodaira surface as an example \cite{Poon1}.

While this development continues to exert its impact on the development of homological aspects of
mirror symmetry, see e.g. 
\cite{CLP} \cite{COP} \cite{Costello-Li} \cite{KKP1} \cite{KKP2} \cite{LiSi}, its degree-2 component
\[
  \sum_{p+q=2}H^q(M, \wedge^p\Theta)\subset  H^\bullet(M, \wedge^\bullet\Theta)
  \]
 of extended
deformation could be  interpreted geometrically by the deformation theory of generalized geometry 
 \cite{Marco} \cite{Hitchin-Generalized CY} \cite{Hitchin-Instantons}. 
 In particular, $H^1(M, \Theta)$ corresponds exactly to the linearized "virtual moduli" space of classical
 deformation theory \cite{Kodaira}. If one finds a de-Rham closed representative in 
 $H^2(M, \wedge^0\Theta)=H^2(M, \cal{O})$, it is due to a $B$-field transformation within 
 the realm of generalized geometry \cite{Marco}. 
 As we will explain in subsequent computation,  closed 2-forms are both $\dbar$-closed and they commute with every element in the 
 Schouten algebra that controls extended deformations, they don't contribute any non-trivial terms in the homological algebra aspect
 of deformation theory. $B$-field transformation is, therefore, a geometrical realization of this algebraic fact. 

However, if an element in $H^2(M, \cal{O})$ is not de-Rham closed, it represents a deformation in the sense of 
 "gerbe" as seen in \cite{Hitchin-Gerbe} and \cite{Hitchin-bundles}.  
 On the other hand, when an element in $H^0(M, \wedge^2\Theta)$ represents what is known as holomorphic
 Poisson structure, it is \it integrable \rm in the sense of extended deformation theory. 
 
 While classical complex deformation has been a century-old subject, holomorphic Poisson structures are also 
 known for some time \cite{B-Marci} \cite{Xu}
 \cite{Polish}. Its recent role in generalized geometry and related deformation theory
 arouse further interests in this subject \cite{FM} \cite{Goto0} \cite{Goto}
 \cite{Hitchin-holomorphic Poisson} \cite{Ran1} \cite{Ran2}. 
 Since generalized geometry encompasses both 
 complex structures and symplectic structures, holomorphic Poisson structure 
 has been particularly useful as a tool to capture both features for 
  holomorphic symplectic
 manifolds \cite{BDV} \cite{Hitchin-Instantons} \cite{Ran1} \cite{Ran2}. 
 In addition, it is now known that holomorphic Poisson structures play a fundamental role in generalized geometry
 \cite{Bailey}. 
 Therefore, there has been work towards 
 understanding holomorphic Poisson cohomology \cite{CGP} \cite{Hong} \cite{Hong-Xu} \cite{Poon2}. 
 
 It is not hard to recognize that holomorphic Poisson cohomology is the hypercomplex of a bi-complex. 
Therefore, theoretically speaking it could be computed by spectral sequence method once a filtration is chosen. 
The first page of the spectral sequences of one of the two natural filtrations of this bi-complex leads to a holomorphic 
version of Lichnerowicz-Poisson cohomology \cite{Lich} \cite{Vaisman}.  The first page of the spectral sequence of another 
natural filtrations leads to the Dolbeault cohomology of polyvector fields $H^q(M, \wedge^p\Theta)$. The $d_1$-map of the 
spectral sequence is simply the restriction of the adjoint action of the holomorphic Poisson structure $\Lambda$, where 
$\adL (-)=\lbra{\Lambda}{-}$. Much of the recent computation of holomorphic Poisson cohomologies takes
the second approach. However as seen in \cite{CGP} \cite{Hong-Xu}
\cite{Xu}  \cite{Poon2}, this spectral sequence rarely degenerates on the second page on manifolds such 
as rational surfaces or Hopf manifolds, let alone
degeneracy on the first page. Therefore, in general 
the spectral sequence approach provides a theoretical but not an effective computational framework.  

On the other hand, nilmanifolds are also known to be rich source of generalized geometry of various kinds \cite{CG}. 
It has long been known that the DeRham cohomology of nilmanifolds is given by algebraic objects \cite{Nomizu}. 
Through a series of work by many, e.g. \cite{CGP}  \cite{Console} \cite{Console-Fino}  \cite{Rolle}, 
we  learn that the Dolbeault cohomology with structure sheaf or tangent sheaf
are  given by algebraic objects  when the holomorphic
Poisson structure is also an algebraic object on a finite-dimensional space. 
The significance of such observation is
that the computation of spectral sequence becomes algebraic as well, and much of the computation becomes
tractable in an elementary manner. 

Based on such foundation, the author and his collaborators recently identify when the spectral sequence degenerates on 
the first page if the complex structure is abelian and the underlying nilmanifold is 2-steps \cite{PS1} \cite{PS2}. It enables an effective 
computation of holomorphic Poisson cohomology. 
They go on to identify the necessary and sufficient condition for the holomorphic Poisson cohomology, as a Gerstenhaber algebra, 
to be isomorphism to case when the Poisson structure is equal to zero. If one considers $t\Lambda$ to be a deformation with deformation 
variable $t\in [0,1]$, it means that the holomorphic Poisson cohomology is invariant under special condition. 
They also demonstrate that the primary Kodaira surface as seen in \cite{Poon1} is an example. 

It takes us back to an observation in \cite{Poon1}. In \cite{Poon1}, the author identifies the Frobenius structure of primary Kodaria 
surfaces. An intriguing observation, after the computation, is that the restriction of the Frobenius structure on the even part of the supermanifold
is trivial. The goal of this paper is to prove that the Frobenius structure restricted to the moduli of a generalized deformation of 
a primary Kodaira manifold of any dimension is trivial. It establishes a conjecture that the same is true for a large class of 2-step nilmanifolds with 
abelian complex structures. We will provide a rigorous statement of our theorem in the next section after we acquire sufficient 
preparation to convey the technical details. In the last section, we will present a conjecture to indicate why and potentially how we may extend
the results in this paper beyond Kodaira manifolds. 

 \section{Deformation Theory}

For any complex manifold $M$, let $\Theta$ be the holomorphic tangent bundle and $\Omega$ the bundle of $(1,0)$-forms. 
We use $\overline\Theta$ and $\Oomega$ to denote their conjugate bundles.
It is well known since Kodaira-Spencer's deformation theory \cite{Kodaira-Spencer} that a deformation of complex
structure is given by $\Gamma$ with 
$\Gamma\in C^\infty(M, \Theta\otimes \Oomega)$ satisfying the 
Maurer-Cartan equation:
\[
\dbar\Gamma+\frac12\lbra{\Gamma}{\Gamma}=0,
\]
where $\lbra{-}{-}$ is the Schouten-Nijenhuis bracket
 \cite{Kodaira-Spencer} 
\cite{Mac}. 
The Schouten-Nijenjuis bracket, or simply the Schouten bracket, could be extended to the
space of sections of the bundle of exterior algebra 
\[
\wedge^\bullet\left(\Theta\oplus \Oomega\right)
=\oplus_{p,q} \left( \wedge^p\Theta \otimes \wedge^q\Oomega \right).
\]
We adopt the notations that $\Theta^p= \wedge^p\Theta$ and
$\Oomega^q= \wedge^q\Oomega$. For a section in $\Theta^p\otimes\Oomega^q$, 
we refer $(p,q)$ as its bi-degree, and $p+q$ its (total) degree. 
We observe that if $U$ and $V$ are $(1,0)$-vector fields, $\overline\alpha$ and 
$\overline\beta$ are $(0,1)$-forms, then $\lbra{U}{V}$ is simply the 
usual Lie bracket between vector fields; $\lbra{U}{\beta}={\cal{L}}_U\beta=d\iota_U\beta+\iota_U d\beta$ is the 
Lie derivative of the $1$-form $\beta$ with respect to the vector field
$U$; and $\lbra{\alpha}{\beta}=0$. It extends to the full exterior algebra
through the rules of Schouten algebra \cite{Mac} \cite{Poon1}. 

On the other hand, $\dbar\alpha$ is the $(0,2)$-part of  $d\alpha$, the exterior differential 
of $\alpha$. $\dbar U$ is the Cauchy-Riemann differentiation as noted in 
\cite{Gau} \cite{Poon1}. The $\dbar$-operator extends to the full exterior algebra
$\wedge^\bullet\left(\Theta\oplus \Oomega\right)$ so that it becomes an exterior differential 
algebra. Due to the integrability of the complex structure, $\dbar\circ\dbar=0$. 
In particular, we will consider the cohomology spaces $H^\bullet=\oplus_kH^k(M)$, where
the $k$-th cohomology space $H^k(M)$ has a bi-degree decomposition. 
\begin{equation}
H^k(M)=\oplus_{p+q=k}H^q(M, \Theta^p). 
\end{equation}

The exterior product
$\wedge$ and the Schouten bracket on 
the space $C^\infty \left( M, \wedge^\bullet(\Theta\oplus \Oomega) \right)$ 
are compatiable in the sense that they form a Gerstenhaber algebra \cite{Gerstenhaber}.
With the $\dbar$-operator, 
\[
\left( C^\infty \left( M, \wedge^\bullet(\Theta\oplus \Oomega) \right), \quad \lbra{-}{-}, \quad
\wedge, \quad \dbar \right)
\]
form a differential Gerstenhaber algebra $DGA(0)$ \cite{Merkulov}
\cite{Poon1}. In particular, for any sections $K$ and $L$ with degree $k$ and $\ell$ respectively, 
\begin{eqnarray}
&\dbar(K\wedge L)=(\dbar K) \wedge L+(-1)^kK\wedge(\dbar L),&\label{dbar wedge}\\
&\dbar\lbra{K}{L}= \lbra{\dbar K}{L}+(-1)^{k+1}\lbra{K}{\dbar L}.\label{dbar bracket}&
\end{eqnarray}
An element 
$\Gamma\in C^\infty\left( M, \wedge^\bullet(\Theta\oplus \Oomega) \right)$ is an \it integrable \rm
 extended  deformation if it satisfies the Maurer-Cartan equation
\begin{equation}\label{MC}
\dbar\Gamma+\frac12\lbra{\Gamma}{\Gamma}=0.
\end{equation}
For such $\Gamma$, the operator 
\begin{equation}
\dbar_\Gamma=\dbar+\lbra{\Gamma}{-}
\end{equation}
satisfies $\dbar_\Gamma\circ \dbar_\Gamma=0$ and hence determines a cohomology
$H_\Gamma^\bullet(M)$ \cite{LWX}. 
Due to (\ref{dbar wedge}) and (\ref{dbar bracket}) above, the cohomology space inherits the structure of a Gerstenhaber 
algebra. 

Roughly speaking, if the collection of integrable extended deformations $\Gamma$ 
forms a "moduli" space when equivalence is accounted for, 
the Frobenius structure on the moduli space measures how the cohomology 
$H_\Gamma^\bullet(M)$ varies when $\Gamma$ varies. 

Denote the linear part  of $\Gamma$ by $\Gamma_1$.   
The linear portion of the Maurer-Cartan equation (\ref{MC}) shows that $\dbar\Gamma_1=0$. In particular, 
it is an element in $H^\bullet(M)$. Conversely, given $\Gamma_1$ in $H^\bullet(M)$, one may attempt the 
Kuranishi's recursive method  to identify an infinite series of obstruction to the existence
of $\Gamma$ \cite{Kuranishi}. In this notes, we work on primary Kodaira manifolds, which will be defined in the next section. 
Instead of working with the full extended deformation, we focus on the portion leading to
 generalized deformations, namely
\[
\Gamma_1\in \oplus_{p+q=2}H^q(M, \Theta^p).
\]
We consider such $\Gamma_1$ as infinitesimal generalized deformation. 
Our goal in this paper is to derive the following theorem. 

\begin{theorem}\label{main} Let $M$ be a primary Kodaira manifold. Every infinitesimal generalized deformation sufficiently close to zero is integrable. 
In other words, for each
 $\Gamma_1\in \oplus_{p+q=2}H^q(M, \Theta^p)$ near $0$,  there exists $\Gamma \in C^\infty\left( M, \wedge^\bullet(\Theta\oplus \Oomega) \right)$ with $\Gamma_1$ being its linear part,
satisfying the Maurer-Cartan equation. Moreover, there exists a natural isomorphism of Gerstenhaber algebras: 
 \[
 \left( H^\bullet_\Gamma(M), \quad \lbra{-}{-}, \quad \wedge \right) 
\quad  \cong \quad
\left( H^\bullet(M), \quad \lbra{-}{-}, \quad \wedge \right).
\]
\end{theorem}

In Theorem \ref{Main} below, we will provide further details regarding how close to zero is sufficient. 
This theorem means that Frobenius structure restricted to the degree-2 portion of the $\ZZ$-grading 
of the moduli space associated to a primary Kodaira manifold is trivial. It generalizes a computation 
on Kodaira surface in \cite{Poon1} from complex 2-dimension to all 
complex dimensions. It also generalizes an observation in \cite{PS2} that the holomorphic Poisson
cohomology of any holomorphic Poisson 
structures on primary Kodaira manifolds of all dimensions, as Gerstenhaber algebra, are isomorphic to 
$H^\bullet(M)$.

\section{Kodaira Manifolds}

We consider primary Kodaira manifolds as nilmanifolds, meaning that they are co-compact quotient of  a simply-connected
nilpotent Lie group. 

Let $\{X_1, Y_1, \dots, X_n, Y_n, Z_1, Z_2\}$ be the basis of a real vector space of 
$2n+2$-dimension. Define on it a Lie bracket $\lbra{-}{-}$ such that the only non-trivial
ones are given by
\begin{equation}
\lbra{X_j}{Y_j}=-\lbra{Y_j}{X_j}=Z_1.
\end{equation}
The Lie bracket turns the concerned vector space in a a Lie algebra, which we denote by $\lieg$. 
Its center $\liec$ is spanned by $Z_1, Z_2$.
Let $G$ be its corresponding simply connected Lie group. Its co-compact quotient is a (primary) Kodaira manifold, which we denote by $M$.
A detailed description of the related lattice and group structure could be found in \cite{GMPP}. 

Consider the linear map $J:\lieg \to \lieg$ such that
\begin{equation}\label{Lie structure}
JX_j=Y_j, \quad JY_j=-X_j, \quad JZ_1=Z_2, \quad JZ_2=-Z_1.
\end{equation}
The operator $J$ becomes an integrable complex structure on $\lieg$. Through left translation, 
it becomes an integrable complex structure on $G$, and descends
on the co-compact quotient which is defined by right action.
For each $1\leq j\leq n$, define $T_j=\frac12 (X_j-i Y_j)$ and $W=\frac12(Z_1-i Z_2)$. 
The space $\lieg^{1,0}=\{T_1, \dots, T_n, W\}$ form the space of $(1,0)$-vectors.
Let $\liet^{1,0}$ and $\liec^{1,0}$ respectively be spanned by $\{T_1, \dots, T_n\}$ and $W$.
Their conjugate spaces are denoted by $\overline{\lieg^{1,0}}$, $\overline{\liet^{1,0}}$ and $\overline{\liec^{1,0}}$ respectively.
Given the way the complex structure $J$ is defined, $\lieg^{1,0}$ is an abelian 
complex algebra while on $\lieg^{1,0}\oplus \overline{\lieg^{1,0}}$
in terms of the basis $\{T_1, \dots, T_n, W, {\overline T}_1. \dots, {\overline T}_n, {\overline W}\}$, 
the structure equations are given by
\begin{equation}
\lbra{T_j}{\overline{T}_j}=\frac14\lbra{X_j-iY_j}{X_j+iY_j}=\frac{i}2\lbra{X_j}{Y_j}=\frac{i}2Z_1=\frac{i}2(W+\overline{W}).
\end{equation}

Let $\{\omega^1, \dots, \omega^n, \rho\}$ be the dual basis with respective to $\{T_1, \dots, T_n, W\}$. 
Denote the complex linear span of its complex conjugation by $\lieg^{0,1}$.
We will also denote the $p$-th exterior product of 
$\lieg^{1,0}$ by $\lieg^{p,0}$,  the $q$-th exterior product of
$\lieg^{0,1}$ by $\lieg^{0,q}$, and denote $\lieg^{p,0}\otimes \lieg^{0,q}$ by $\lieg^{p,q}$. Then
\begin{equation}
\wedge^k(\lieg^{1,0}\oplus {\overline{\lieg^{1,0}}})=\sum_{p+q=k}\lieg^{p,q}.
\end{equation}
$\liet^{p,q}$ and $\liec^{p,q}$ are similarly defined.
In this setting the Lie algebra differential $\dbar: \lieg^{0,1}\to \lieg^{0,2}$ is determined by
\begin{equation}
d\oomega^j=0 \quad \mbox{ for all } j, \quad \mbox{ and } \quad d\orho=-\frac{i}2 \sum_{j=1}^n\omega^j\wedge\oomega^j.
\end{equation}
Apparently, $d\rho=d\orho$ and 
\begin{equation}
\dbar\omega^j=0 \quad  \mbox{ for all } j, \quad \mbox{ and } \quad \dbar\orho=0.
\end{equation}
To compute $\dbar: \lieg^{1,0} \to \lieg^{1,1}$, we note that the Cauchy-Riemann operator is given by
\[
\dbar_{\overline B}A=\lbra{{\overline B}}{A}^{1,0}
\]
for any $A\in \lieg^{1,0}$ and ${\overline B}\in {\overline{\lieg^{1,0}}}$ \cite{Gau}. Therefore for all $T\in \lieg^{1,0}$,
\begin{equation}
\dbar T=\sum_{j=1}^n\lbra{{\overline T}_j}{T}^{1,0}\wedge \oomega^j.
\end{equation}
It follows that
\begin{equation}
\dbar W=0, \quad \mbox{ and for all }  j,  \quad  \dbar T_j=-\frac{i}2 W\wedge \oomega^j.
\end{equation}

Through contraction $\iota_{T_j}d\rho=-\frac{i}2\oomega^j$,
the maps $d\rho$ and $d\orho$ are treated as linear maps.
\begin{equation}
d\rho,  d\orho \quad : \liet^{1,0} \to \liet^{0,1}.
\end{equation}
As explained in \cite{Poon1}, 
\begin{equation}
\Ger=\left( \oplus_{p,q}\lieg^{p,q}, \quad \lbra{-}{-}, \quad \wedge, \quad \dbar \right) 
\end{equation} 
is a finite-dimensional Gerstenbaber sub-algebra contained in 
\[
\left( C^\infty(M, \wedge^\bullet (\Theta\oplus\Oomega) ),   \quad \lbra{-}{-}, \quad \wedge,  \quad \dbar   \right).
\]
Moreover, the inclusion map is a quasi-isomorphism, meaning that it induces an isomorphism on 
cohomology level \cite[Proposition 6]{Poon1}. Therefore, from now on we focus our computation on the finite-dimensional 
algebra $\Ger$. 

The vector space $\oplus_{p,q}\lieg^{p,q}$ with the Schouten bracket form a \it Schouten algebra. \rm
\[
\Sch=\left(  \oplus_{p,q}\lieg^{p,q}, \quad \lbra{-}{-} \right).
\]
The structure of $\Sch$  is generated by the bracket restricted to degree-1 elements
\begin{equation}
\Sch^1=\lieg^{1,0}\oplus \lieg^{0,1},
\end{equation}
where $\lbra{T}{\oom}$ by definition is equal to $\iota_Td\oom$. As the algebra $\lieg^{1,0}$ 
is abelian and all elements in $\liet^{0,1}$ is $d$-closed, the
only non-trivial bracket is due to
\begin{equation}\label{non trivial bracket}
\lbra{T_j}{\orho}=\iota_{T_j}d\orho=-\frac{i}2 \oom^j.
\end{equation}
Referencing to the structure of the Schouten algebra $\Sch$, we note the following.
\begin{lemma}\label{degree 1 center}
The spaces $\liec^{1,0}$ and $\liet^{0,1}$  are in the center of the Schouten algebra.
\end{lemma}

In subsequent analysis, we focus on the even part of the Schouten algebra,
\begin{equation}
\Sch^{\mbox{even}}=\oplus_{p+q=\mbox{even}}\lieg^{p,q}.
\end{equation}
We start with $\Sch^2=\oplus_{p+q=2}\lieg^{p,q}=\lieg^{2,0}\oplus \lieg^{1,1}\oplus \lieg^{0,2}.$
It contains $\liet^{1,1}$. We further define the following for all $1\leq i, j\leq n$.
\begin{equation}\label{def phi psi}
\phi_{ij}=\frac12 (T_i\wedge \oom^j+T_j\wedge \oom^i),
\quad
\psi_{ij}=\frac12 (T_i\wedge \oom^j-T_j\wedge \oom^i).
\end{equation}
Apparently, $\phi_{ij}=\phi_{ji}$ and $\psi_{ij}=-\psi_{ji}$ for all $i, j$. We denote the respective linear spans by
$\odot^{1,1}$ and $\triangle^{1,1}$ so that
\begin{equation}
\liet^{1,1}=\odot^{1,1}\oplus \triangle^{1,1}.
\end{equation}
So we have
\begin{itemize}
\item $\lieg^{2,0}=(\liec^{1,0}\otimes \liet^{1,0}) \oplus \liet^{2,0}.$
\item $\lieg^{1,1}=\liec^{1,1}\oplus (\liec^{1,0}\otimes \liet^{0,1} )\oplus  (\liet^{1,0}\otimes \liec^{0,1})\oplus \odot^{1,1}\oplus \triangle^{1,1}.$
\item $\lieg^{0,2}=(\liec^{0,1}\otimes \liet^{0,1} )\oplus \liet^{0,2}. $
\end{itemize}

Since $\liec^{1,0}$ and $\liet^{0,1}$ are in the center of the Schouten algebra, 
 as far as Schouten bracket is concerned, 
\begin{itemize}
\item $\lieg^{2,0}\equiv_s(\liec^{1,0}\otimes \liet^{1,0}) \oplus \liet^{2,0}.$
\item $\lieg^{1,1}\equiv_s \liec^{1,1}\oplus  (\liet^{1,0}\otimes \liec^{0,1})\oplus \odot^{1,1}\oplus \triangle^{1,1}.$
\item $\lieg^{0,2}\equiv_s\liec^{0,1}\otimes \liet^{0,1} . $
\end{itemize}

In addition, due to (\ref{non trivial bracket}) and definitions in (\ref{def phi psi})
\begin{equation}
\lbra{T_j\wedge T_k}{\orho} = T_j\wedge\lbra{T_k}{\orho}-T_k\wedge \lbra{T_j}{\orho}
=-\frac{i}2(T_j\wedge\oom^k-T_k\wedge\oom^j)=-i\psi_{jk}.
\label{tjk orho}
\end{equation}

Below we compute the Schouten bracket among terms in $\Sch^2$. Firstly, 
\begin{equation}
\lbra{\lieg^{2,0}}{\lieg^{2,0}}=0
\end{equation}
because the complex structure is abelian.
Since $\liec^{1,0}$ is in the center and is one-dimensional,
\begin{equation}
\lbra{\liec^{1,0}\otimes \liet^{1,0}}{\liec^{1,1}}=0.
\end{equation}
To compute $\lbra{\liec^{1,0}\otimes \liet^{1,0}}{\liet^{1,0}\otimes \liec^{0,1}}$, consider $W\wedge T_j$ 
in $\liec^{1,0}\otimes \liet^{1,0}$ and
$T_k\wedge \orho$ in $\liet^{1,0}\otimes \liec^{0,1}$.
\begin{equation}
\lbra{W\wedge T_j}{T_k\wedge \orho}=-W\wedge \lbra{T_j}{ \orho}\wedge T_k=-\frac{i}2W\wedge T_k\wedge \oom^j.
\end{equation}
In particular,
\begin{equation}
\lbra{\liec^{1,0}\otimes \liet^{1,0}}{\liet^{1,0}\otimes \liec^{0,1}}\subseteq \liec^{1,0}\otimes \liet^{1,1}.
\end{equation}
Given (\ref{non trivial bracket}),
\begin{equation}
\lbra{\liec^{1,0}\otimes \liet^{1,0}}{\odot^{1,1}}=0
\quad  \mbox{ and }  \quad
\lbra{\liec^{1,0}\otimes \liet^{1,0}}{\triangle^{1,1}}=0.
\end{equation}
i.e.
\begin{equation}\label{c10t10t11}
\lbra{\liec^{1,0}\otimes \liet^{1,0}}{\liet^{1,1}}=0.
\end{equation}

Regarding $\lbra{\liec^{1,0}\otimes \liet^{1,0}}{\liec^{0,1}\otimes\liet^{0,1}}\subseteq \liec^{1,0}\otimes\liet^{0,2}$, we have
\begin{equation}\label{wtj orho oomk}
\lbra{W\wedge T_j}{\orho\wedge\oomega^k} =
W\wedge \lbra{T_j}{\orho}\wedge\oomega^k=-\frac{i}2 W\wedge \oom^j\wedge\oomega^k.
\end{equation}
In summary,
\begin{equation}
\lbra{\liec^{1,0}\otimes \liet^{1,0}}{\lieg^{1,1}}
=\lbra{\liec^{1,0}\otimes \liet^{1,0}}{\liec^{1,1}\oplus  (\liet^{1,0}\otimes \liec^{0,1})}
\subseteq (\liec^{1,0}\otimes \liet^{1,1}) \oplus (\liec^{1,0}\otimes\liet^{0,2}).
\end{equation}

Similarly,
\begin{equation}
\lbra{\liet^{2,0}}{\odot^{1,1}}=0, \quad  \lbra{\liet^{2,0}}{\triangle^{1,1}}=0, \quad
\lbra{\liet^{2,0}}{\liec^{0,1}\otimes\liet^{0,1}}=\liet^{0,1}\otimes {\triangle}^{1,1},
\end{equation}
because by (\ref{tjk orho}),
\begin{eqnarray}
\lbra{T_j\wedge T_k}{W\wedge \orho}
&=&i W\wedge {\psi}_{jk},
\\
\lbra{T_j\wedge T_k}{T_\ell\wedge \orho}
&=&iT_\ell\wedge {\psi}_{jk},
\\
\lbra{T_j\wedge T_k}{ \orho \wedge \oom^\ell}
&=&-i{\psi}_{jk}\wedge \oom^\ell.
\end{eqnarray}
Next $\lbra{\liec^{1,1}}{\liec^{1,1}}=0$. 
Since
\begin{eqnarray}
\lbra{W\wedge\orho}{T_j\wedge\orho} &=&W\wedge\lbra{\orho}{T_j}\wedge\orho
=-\frac{i}2 W\wedge\orho\wedge \oom^j,\\
\lbra{\orho}{\phi_{jk}}
&=&
\frac12
\left(
\lbra{\orho}{T_j}\wedge \oom^k+\lbra{\orho}{T_k}\wedge \oom^j
\right)
=0,
\label{orho phijk}
\\
\lbra{\orho}{\psi_{jk}}
&=&
\frac{i}4
\left(
\oom^j\wedge \oom^k-\oom^k\wedge \oom^j
\right)
=\frac{i}2 \oom^j\wedge \oom^k
\label{orho psijk}
\end{eqnarray}
for all $i, j$, we observe that 
\[
\lbra{\liec^{1,1}}{\liet^{1,0}\otimes \liec^{0,1}} = \liec^{1,1}\otimes \liet^{0,1} ,
\quad
\lbra{\liec^{1,1}}{\odot^{1,1}}=0,
\quad
\lbra{\liec^{1,1}}{\triangle^{1,1}}=\liec^{1,0}\otimes \liet^{0,2}.
\]
As a consequence, we note the following
\begin{lemma}\label{degree 2 center} The spaces $\odot^{1,1}$, $\liec^{1,0}\otimes \liet^{0,1}$ and $\liet^{0,2}$ in $\Sch^2$ are in the center of the Schouten algebra
$\Sch$, and hence up to equivalence in Schouten algebra,
\begin{itemize}
\item $\lieg^{2,0}\equiv_s(\liec^{1,0}\otimes \liet^{1,0}) \oplus \liet^{2,0}.$
\item $\lieg^{1,1}\equiv_s \liec^{1,1}\oplus  (\liet^{1,0}\otimes \liec^{0,1})\oplus \triangle^{1,1}.$
\item $\lieg^{0,2}\equiv_s\liec^{0,1}\otimes \liet^{0,1} . $
\end{itemize}
\end{lemma}
It follows (\ref{orho psijk}) that
\begin{equation}
\lbra{W\wedge\orho}{\psi_{jk}}
=
\frac{i}2 W \wedge \oom^j\wedge \oom^k
\label{c11tri11}
\end{equation}
for all $i, j$.
Also,
$\lbra{\liec^{1,1}}{\liec^{0,1}\otimes\liet^{0,1}}=0.$
By (\ref{non trivial bracket}),
\begin{eqnarray}
\lbra{T_j\wedge\orho}{T_k\wedge\orho}
&=&
T_j\wedge\lbra{\orho}{T_k}\wedge\orho+\orho\wedge\lbra{T_j}{\orho}\wedge T_k
\nonumber\\
&=&
\frac{i}2
\left( T_j\wedge \oom^k\wedge\orho-\orho\wedge \oom^j\wedge T_k
\right)
=i{\psi}_{jk}\wedge\orho.
\end{eqnarray}
By  (\ref{orho psijk}),
\begin{equation}
\lbra{T_\ell\wedge\orho}{{\psi}_{jk}}
=\frac{i}2 T_\ell\wedge \oom^j\wedge\oom^k.
\end{equation}
It follows that 
\begin{eqnarray*}
\lbra{\liet^{1,0}\otimes \liec^{0,1}}{\liet^{1,0}\otimes \liec^{0,1}}&=&{\triangle}^{1,1}\otimes \liec^{0,1},
\\
\lbra{\liet^{1,0}\otimes \liec^{0,1}}{\triangle^{1,1}}&=& \liet^{1,0}\otimes\tri^{1,1},
\\
\lbra{\liet^{1,0}\otimes \liec^{0,1}}{\liec^{0,1}\otimes\liet^{0,1}}&=&\liec^{0,1}\otimes\liet^{0,2}.
\end{eqnarray*}
Next, it is clear that
$ \lbra{\triangle^{1,1}}{\triangle^{1,1}}=0.$
We also note that
$\lbra{\triangle^{1,1}}{\liec^{0,1}\otimes \liet^{0,1}}=\liet^{0, 3}$
because
\begin{equation}\label{tri11c01t01}
\lbra{\psi_{jk}}{\orho\wedge \oom^\ell}
=
\frac{i}2
\oom^j\wedge \oom^k\wedge \oom^\ell.
\end{equation}
Finally,
$
\lbra{\liec^{0,1}\otimes\liet^{0,1}}{\liec^{0,1}\otimes\liet^{0,1}}=0.$
Therefore, we could now summarize the Schouten bracket on $\Sch^2$ in a $5\times 5$-symmetric matrix.
\begin{proposition}\label{degree 2 sch}
All non-zero terms of Schouten bracket among elements in $\Sch^2$ are given below.
\[
\begin{array}{|c||c|c|c|c|c|c|}
\hline
 & \liec^{1,0}\otimes \liet^{1,0}  & \liet^{2, 0}  & \liec^{1,1}  & \liet^{1,0}\otimes \liec^{0,1}  & {\triangle}^{1,1}  & \liec^{0,1}\otimes \liet^{0,1} \\
     \hline \hline
\liec^{1,0}\otimes \liet^{1,0} & 0 & 0 & 0 & \liec^{1,0}\otimes\liet^{1,1} &  0 & \liec^{1,0}\otimes\liet^{0,2}
\\
\hline
 \liet^{2, 0} & 0 & 0 & \liec^{1,0}\otimes{\triangle}^{1,1} & \liet^{1,0}\otimes{\triangle}^{1,1}   & 0 & \liet^{0,1}\otimes{\triangle}^{1,1}\\
 \hline
 \liec^{1,1} & 0 & \liec^{1,0}\otimes {\triangle}^{1,1} &0 & \liec^{1,1}\otimes\liet^{0,1} & \liec^{1,0}\otimes \liet^{0,2} & 0\\
 \hline
 \liet^{1,0}\otimes\liec^{0,1} & \liec^{1,0}\otimes\liet^{1,1} & \liet^{1,0}\otimes{\triangle}^{1,1} &\liec^{1,1}\otimes\liet^{0,1} &{\triangle}^{1,1}\otimes\liec^{0,1}
  & \liet^{1,2} & \liec^{0,1}\otimes\liet^{0,2}\\
 \hline
 \triangle^{1,1} & 0 & 0 & \liec^{1,0}\otimes\liet^{0,2} & \liet^{1,2}  & 0 & \liet^{0,3}\\
 \hline
 \liec^{0,1}\otimes\liet^{0,1} & \liec^{1,0}\otimes\liet^{0,2} & \liet^{0,1}\otimes {\triangle}^{1,1} & 0 & \liec^{0,1}\otimes \liet^{0, 2}  & \liet^{0, 3} & 0\\
 \hline
\end{array}
\]
\begin{center}{ Table 1: Schouten bracket among degree-2 elements}
\end{center}
\end{proposition}

\section{Exterior differential algebra}
In this section, we consider the exterior differential algebra structure within the Gerstenhaber algebra:
$
\left( 
 \oplus_{p,q}\lieg^{p,q},   \wedge, \dbar 
\right).$
All are dictated by the equations
\begin{equation}\label{full structure}
\dbar W=0, \quad \dbar T_j=-\frac{i}2 W\wedge \oom^j, \quad \dbar\orho =0, \quad \dbar \oom^j=0
\end{equation}
for all $1 \leq j\leq n$.
We now begin to compute the cohomology.
\begin{equation}
H^{p,q}=\frac{\ker \dbar : \lieg^{p,q} \to \lieg^{p, q+1} }{{\image } \ \dbar: \lieg^{p,q-1} \to \lieg^{p, q}}
\end{equation}

To facilitate further computation in deformation theory, we consider a vector space decomposition $\lieg^{p,q}$ in three components.
Its so-called harmonic part  is isomorphism to $H^{p, q}$. Its $\dbar$-exact part is $D^{p, q}=\image \dbar: \lieg^{p, q-1} \to \lieg^{p,q}$.
Taking advantage of the finite-dimensional situation, 
we will choose a base so that the remaining part is denoted by $G^{p, q}$. 
In respect to classical deformation theory, we address $G^{p, q}$ the Green's part. We consider
\[
\lieg^{p,q}=H^{p,q}\oplus D^{p, q}\oplus G^{p,q}
\]
a Hodge decomposition of the type-$(p,q)$ space $\lieg^{p,q}$. 
From observations above (\ref{full structure}), we conclude the following.
\begin{lemma}\label{type 10}\label{type 01}
The Hodge decomposition of degree-1 cohomology is given below.
\begin{eqnarray}
&H^{1,0}=\liec^{1,0}, \quad D^{1,0}=0, \quad G^{1,0}=\liet^{1,0},&\\
&H^{0,1}=\liec^{0,1}\oplus \liet^{0,1}, \quad D^{0,1}=0, \quad G^{0,1}=0.&
\end{eqnarray}
Moreover, the map $\dbar: \liet^{1,0}\to \liec^{1,0}\otimes\liet^{0,1}$ where
\begin{equation}
 \dbar T_j=-\frac{i}2 W\wedge \oom^j,
\end{equation}
is an isomorphism.
\end{lemma}

Due to the computation above, it is now apparent that for all $T$, $\dbar(W\wedge T)=0$. However,
\begin{eqnarray*}
\dbar (T_j\wedge T_k) &=&-\frac{i}2 \left( W\wedge \oom^j\wedge T_k-T_j\wedge W\wedge \oom^k\right)\\
&=& -\frac{i}2 W\wedge (T_j\wedge \oom^k -T_k\wedge \oom^j)=-i W\wedge \psi_{jk}.
\end{eqnarray*}

\begin{lemma}\label{type 20} The Hodge decomposition of Type-(2,0) space is given below.
\begin{equation}
H^{2,0}=\liec^{1,0}\otimes\liet^{1,0}, \quad D^{2,0}=0, \quad G^{2,0}=\liet^{2,0}.
\end{equation}
Moreover, the map $\dbar: \liet^{2,0} \to \liec^{1,0}\otimes \triangle^{1,1}$ where
\begin{equation}
  \dbar (T_j\wedge T_k)=-i W\wedge \psi_{jk}
\end{equation}
is an isomorphism.
\end{lemma}

As vector spaces,
\begin{equation}
\lieg^{1,1}=\liec^{1,1}\oplus (\liec^{1,0}\otimes\liet^{0,1}) \oplus (\liet^{1,0}\otimes\liec^{0,1}) \oplus \odot^{1,1}\oplus \tri^{1,1}.
\end{equation}
It is apparent that $W\wedge\orho$ is $\dbar$-closed but not $\dbar$-exact. Due to Lemma \ref{type 10}, the second summand above is $\dbar$-exact.
The same implies that for all $T_j$ in $\liet^{1,0}$,
\begin{equation}
\dbar (T_j\wedge \orho)=-\frac{i}2 W\wedge \oom^j\wedge \orho.
\end{equation}
Regarding the remaining summands, we observe that
\begin{eqnarray}
\dbar(T_j\wedge \oom^k) &=&-\frac{i}2 W\wedge \oom^j\wedge \oom^k.\\
\dbar \phi_{jk}&=&-\frac{i}4 W\wedge \left(\oom^j\wedge \oom^k+\oom^k\wedge \oom^j \right)=0,\\
\dbar\psi_{jk}&=& -\frac{i}4 W\wedge \left(\oom^j\wedge \oom^k-\oom^k\wedge \oom^j \right)=-\frac{i}2W\wedge \oom^j\wedge \oom^k.
\label{dbar psijk}
\end{eqnarray}

\begin{lemma}\label{type 11} The Hodge decomposition of type-(1,1) space is given below.
\begin{equation}
H^{1,1}=\liec^{1,1}\oplus \odot^{1,1}, \quad
D^{1,1}=\liec^{1,0}\otimes \liet^{0,1}=\dbar\liet^{1,0}, \quad
G^{1,1}=(\liet^{1,0}\otimes \liec^{0,1})\oplus \tri^{1,1}.
\end{equation}
Moreover, the following maps are isomorphisms:
\begin{equation}
\dbar: \liet^{1,0}\otimes \liec^{0,1}\to \liec^{1,1}\otimes\liet^{0,1}, \qquad
\dbar: \tri^{1,1}\to \liec^{1,0}\otimes\liet^{0,2}
\end{equation}
\end{lemma}

On the other hand, the following is clear from the structure equations.
\begin{lemma}\label{type 02}
The Hodge decomposition of type-(0,2) space is given below.
\begin{equation}
H^{0,2}=\lieg^{0,2}=(\liec^{0,1}\otimes\liet^{0,1})\oplus\liet^{0,2},
\quad
D^{0,2}=0, \quad G^{0,2}=0.
\end{equation}
\end{lemma}

Define $H^k=\oplus_{p+q=k}H^{p,q}$, $D^k=\oplus_{p+q=k}D^{p,q}$, and $G^k=\oplus_{p+q=k}G^{p,q}$ respectively, we summarize our observations in the
last few lemmas as below.
\begin{proposition}\label{degree 2 Hodge} The Hodge decomposition of degree-2  space is given below.
\begin{eqnarray*}
H^2&=&(\liec^{1,0}\otimes\liet^{1,0})\oplus \liec^{1,1}\oplus \odot^{1,1} \oplus (\liec^{0,1}\otimes\liet^{0,1}) \oplus\liet^{0,2}\\
D^2&=&\liec^{1,0}\otimes \liet^{0,1}=\dbar\liet^{1,0}\\
G^2&=& \liet^{2,0}\oplus (\liet^{1,0}\otimes \liec^{0,1})\oplus \tri^{1,1}.
\end{eqnarray*}
Moreover,  $ \dbar\liet^{2,0} = \liec^{1,0}\otimes \triangle^{1,1}$,
$ \dbar(\liet^{1,0}\otimes \liec^{0,1})= \liec^{1,1}\otimes\liet^{0,1}$, and
$ \dbar\tri^{1,1}= \liec^{1,0}\otimes\liet^{0,2}.$
\end{proposition}

\section{Solving the Extended Maurer-Cartan Equation}
In view of Proposition \ref{degree 2 sch},  up to non-trivial Schouten brackets 
\[
H^2\equiv_s (\liec^{1,0}\otimes \liet^{1,0})\oplus \liec^{1,1}\oplus (\liec^{0,1}\otimes \liet^{0,1}),
\]
and the 
only non-trivial Schouten bracket among elements in  $H^2$ is due to a single type of  brackets, namely
\[
\lbra{\liec^{1,0}\otimes \liet^{1,0}}{\liec^{0,1}\otimes \liet^{0,1}}\subseteq \liec^{1,0}\otimes \liet^{0,2}=\dbar{\tri}^{1,1}.
\]
To be concrete, by (\ref{wtj orho oomk}) and (\ref{dbar psijk}),
\begin{equation}\label{dbar psijk}
 \lbra{W\wedge T_j}{\orho\wedge \oom^k}=
 -\frac{i}2 W\wedge \oom^j\wedge \oom^k=\dbar\psi_{jk}.
 \end{equation}

As elements in $\tri^{1,1}$ appear as potentials of Schouten bracket between harmonic elements in $\Sch^2$, it is necessary to compute the Schouten bracket
between elements in $H^2$ and $\tri^{1,1}$ as well as Schouten bracket between two elements in $\tri^{1,1}$. However, we have seen that the latter is equal to zero.

Computation in (\ref{c10t10t11}) leads to $\lbra{\tri^{1,1}}{\liec^{1,0}\otimes\liet^{1,0}}=0$. 
Identities (\ref{c11tri11}) and (\ref{dbar psijk}) lead to $\lbra{\tri^{1,1}}{\liec^{1,1}}=\liec^{1,0}\otimes \liet^{0,2}=\dbar\tri^{1,1}$. In fact,
\begin{equation}\label{w orho psijk}
\lbra{W\wedge\orho}{\psi_{ij}}=\frac{i}2 W\wedge\oom^j\wedge\oom^k=-\dbar\psi_{jk}.
\end{equation}
Similarly, (\ref{tri11c01t01}) leads to $\lbra{\tri^{1,1}}{\liec^{0,1}\otimes\liet^{0,1}}= \liet^{0,3}$, i.e.
\begin{equation}\label{psijk orho oom}
\lbra{\psi_{jk}}{\orho\wedge \oom^\ell}=\frac{i}2\oom^j\wedge\oom^k\wedge\oom^\ell.
\end{equation}

On the other hand, elements in $\liet^{0,3}$ are in the center of the Schouten algebra, and they are always in the harmonic part.
Therefore, we have the following observation.
\begin{proposition}
Suppose that $\Gamma_1\in H^2$. There exists $\Psi\in \tri^{1,1}$ and $\vec{\partial}\in \liet^{0,3}$ such that
\begin{equation}
\dbar (\Gamma_1+\Psi)+\frac12\lbra{\Gamma_1+\Psi}{\Gamma_1+\Psi}+\vec{\partial}=0.
\end{equation}
\end{proposition}

After Merkulov \cite{Merkulov}, the element $\vec{\partial}$ above is called the Chen vector and 
it constitutes the obstruction for $\Gamma_1$ to be \it integrable \rm in such a way that
$\Gamma=\Gamma_1+\Psi$ solves the  Maurer-Cartan equation.

We use Kuranishi recursive formula to solve the Maurer-Cartan
equation. Let
\[
\Gamma=\sum_{m=1}^{\infty}t^{m}\Gamma_{m}.
\]
Substitute it into the Maurer-Cartan equation, we have
\begin{eqnarray*}
& \overline{\partial}\left(  \sum_{m=1}^{\infty}t^{m}\Gamma_{m}\right)
+\frac{1}{2}\lbra{ \sum_{k}t^{k}\Gamma_{k}}{\sum_{\ell}t^{\ell}\Gamma_{\ell}} \\
& =t^{m}\left(  \overline{\partial}\Gamma_{m}+\frac{1}{2}\sum_{k+\ell
=m}\lbra{ \Gamma_{k}}{\Gamma_{\ell}}  \right)  .
\end{eqnarray*}
Our goal is that for all $m\geq1$
\[
\overline{\partial}\Gamma_{m}+\frac{1}{2}\sum_{k+\ell=m}\lbra{ \Gamma_{k}}{\Gamma_{\ell}}  =0.
\]
Whenever $\sum_{k+\ell=m}\left[  \Gamma_{k},\Gamma_{\ell}\right]  $ has a
component failing to be $\overline{\partial}$-exact, we compensate it by
$\overrightarrow{\partial}_{m}$. We address it as the $m$-th order term of
the Chen vector. Therefore, we always have
\[
\overline{\partial}\Gamma_{m}+\frac{1}{2}\sum_{k+\ell=m}
\lbra{\Gamma_{k}}{\Gamma_{\ell}}  +\overrightarrow{\partial}_{m}=0.
\]
Define
\[
\overrightarrow{\partial}=\sum_{m=2}^{\infty}t^{m}\overrightarrow{\partial
}_{m},
\]
then
\[
\overline{\partial}\Gamma+\frac{1}{2}\lbra{ \Gamma}{\Gamma}
+\overrightarrow{\partial}=0.
\]
The Kuranishi space of solution of the Maurer-Cartan equation is given by
$\ker\overrightarrow{\partial}.$
In particular,
\[
\overline{\partial}\Gamma_{1}=0, \quad
\frac{1}{2}\sum_{k+\ell=2}\lbra{\Gamma_{k}}{\Gamma_{\ell}}  =-\overline{\partial}\Gamma_{2}-\vec{\partial}_2
\]
where $-\vec{\partial}_2$ is equal to the harmonic part of $\frac{1}{2}\sum_{k+\ell=2}\left[  \Gamma_{k},\Gamma_{\ell}\right]$.
Now, set
$T=\sum_{i}\lambda_{i}T_{i}$, $\oom=\sum_{j}\alpha_{j}\oom^j$. Then
\begin{eqnarray}
\Gamma_{1}
&=&\lambda_{j} W\wedge T_{j}
+\gamma W\wedge\overline{\rho}
+\alpha_{j} \overline{\rho}\wedge\oom^j
+\sum_{j\leq k}\gamma_{jk}\phi_{jk}
+\sum_{j<k}
\beta_{jk} \oom^j\wedge \oom^k
\nonumber \\
&=& W\wedge T+\gamma W\wedge\overline{\rho}
+\overline{\rho}\wedge\oom
+\sum_{j\leq k}\gamma_{jk}\phi_{jk}
+\sum_{j< k}\beta_{jk} \oom^j\wedge \oom^k .
\end{eqnarray}
Since both $\phi_{jk}$ and $\oom^j\wedge\oom^k$ are $\dbar$-closed and in the center of the Schouten algebra, in subsequent
calculation of Kuranishi recursive formula, $\Gamma_1$ is equivalent to
\begin{equation}
\Gamma_1\equiv_s W\wedge T+\gamma W\wedge\overline{\rho}
+\overline{\rho}\wedge\oom.
\end{equation}
By Proposition \ref{degree 2 sch} and identity (\ref{dbar psijk}),
\begin{eqnarray*}
&& \frac12\lbra{\Gamma_1}{\Gamma_1}=\lbra{W\wedge T}{\orho\wedge\oom}
\\
&=&\sum_{j,k}\lambda_j\alpha_k W\wedge\lbra{T_j}{\orho}\wedge \oom^k=-\frac{i}2 \sum_{j,k}\lambda_j\alpha_k W\wedge \oom^j\wedge\oom^k\\
&=&
\dbar\left(\sum_{j,k}\lambda_j\alpha_k\psi_{jk}\right).
\end{eqnarray*}
Therefore,
\[
\Gamma_2=-\sum_{j,k=1}^n\lambda_j\alpha_k\psi_{jk} \quad \mbox{ and } \quad {\vec\partial}_2=0.
\]
Next,
\begin{eqnarray*}
&& \frac12\sum_{k+\ell=3}\lbra{\Gamma_k}{\Gamma_\ell}=\frac12\left( \lbra{\Gamma_1}{\Gamma_2}+\lbra{\Gamma_2}{\Gamma_1}\right)
= \lbra{\Gamma_1}{\Gamma_2}\\
&=&-\sum_{j,k=1}^n\lambda_j\alpha_k
\left( \lbra{W\wedge T}{\psi_{jk}} +\gamma \lbra{W\wedge\orho}{\psi_{jk}}+\lbra{\orho\wedge \oom}{\psi_{jk}}\right)\\
&=&-\sum_{j,k=1}^n\lambda_j\alpha_k
\left( \gamma \lbra{W\wedge\orho}{\psi_{jk}}+\sum_{\ell}\alpha_\ell\lbra{ \orho\wedge \oom^\ell}{\psi_{jk}}\right).
\end{eqnarray*}
By (\ref{w orho psijk}) and (\ref{psijk orho oom}), it is equal to
\begin{eqnarray*}
&=&\gamma \sum_{j,k=1}^n\lambda_j\alpha_k \dbar\psi_{jk}
+\frac{i}2\sum_{j,k=1}^n\lambda_j\alpha_k \oom^j\wedge\oom^k\wedge\oom
\\
&=&\gamma \sum_{j,k=1}^n\lambda_j\alpha_k \dbar\psi_{jk}+\frac{i}2\left(\sum_{j,k=1}^n\lambda_j \oom^j\right)\wedge\oom \wedge\oom\\
&=&-\gamma \dbar\Gamma_2.
\end{eqnarray*}
It means that
\begin{equation}
\Gamma_3=\gamma\Gamma_2 \quad \mbox{ and }  \quad {\vec\partial}_3=0.
\end{equation}
Moreover, we obtain
\begin{equation}\label{induction}
\lbra{\Gamma_1}{\Gamma_2}=-\gamma \dbar\Gamma_2.
\end{equation}
\begin{lemma}
For all $k\geq 3$,
$\Gamma_k=\gamma^{k-2}\Gamma_{2}, \quad {\vec\partial}_k=0. $
\end{lemma}
\bproof We have seen that the above statement is true when $k=3$. Fix any $j$ and suppose that this statement holds for all $n$ such that $j\geq n\geq 3$.
To compute $\Gamma_{j+1}$ and $\vec{\partial}_{j+1}$, we consider
\begin{eqnarray*}
&&\frac12\sum_{k+\ell=j+1}\lbra{\Gamma_k}{\Gamma_\ell}
\\
&=&\frac12 \lbra{\Gamma_{j}}{\Gamma_1}+\frac12\sum_{k+\ell=j+1, k,\ell\geq 2}\lbra{\Gamma_k}{\Gamma_\ell}+\frac12\lbra{\Gamma_1}{\Gamma_j}
\\
&=&\lbra{\Gamma_{j}}{\Gamma_1} +\frac12\sum_{k+\ell=j+1, k,\ell\geq 2}\lbra{\gamma^{k-2}\Gamma_2}{\gamma^{\ell-2}\Gamma_2}\\
&=&\lbra{\Gamma_{j}}{\Gamma_1}=\gamma^{j-2}\lbra{\Gamma_{2}}{\Gamma_1}
\end{eqnarray*}
because $\lbra{\tri^{1,1}}{\tri^{1,1}}=0$. In addition, by (\ref{induction}) it is equal to $-\gamma^{j-1}\dbar\Gamma_2$ as claimed.
\eproof

Since ${\vec\partial}_k=0$ for all $k$, $\vec\partial=0$. Moreover, 
\begin{eqnarray}
&& \Gamma(t)=\sum_{m=1}^\infty t^m\Gamma_m
= t\Gamma_1+\sum_{m=2}^\infty t^m \gamma^{m-2}\Gamma_2=t\Gamma_1+\sum_{m=2}^\infty t^m \gamma^{m-2}\Gamma_2
\nonumber\\
&=&t\Gamma_1+t^2(\sum_{m=2}^\infty t^{m-2} \gamma^{m-2})\Gamma_2=t\Gamma_1+\frac{t^2}{1-t\gamma}\Gamma_2
\end{eqnarray}
when $|t\gamma|<1$.

\begin{theorem}\label{Main} Let $\Gamma_1$ in $H^2$ be given by
\[
\Gamma_{1}
=\sum_j\lambda_{j}\left(  W\wedge T_{j}\right)
+\gamma\left(W\wedge\overline{\rho}\right)
+\sum_k\alpha_{k}\left(  \overline{\rho}\wedge\oom^k\right)
+\sum_{j\leq k}\gamma_{jk}\phi_{jk}
+\sum_{j<k}
\beta_{jk}\left(  \oom^j\wedge \oom^k\right).
\]
If $|\gamma|<1$ and $\Gamma_2=-\sum_{j,k=1}^n\lambda_j\alpha_k\psi_{jk}$, 
\[
\Gamma=\Gamma_1+\frac{1}{1-\gamma}\Gamma_2
\]
is a solution of the Maurer-Cartan equation.
\end{theorem}
It recovers an observation in \cite[Theorem 12]{Poon1}, and constitutes the foundation for our main observation next. 

\section{Frobenius Structure}
We begin to examine the Frobenius structure associated to the primary Kodaira manifolds. Once again, 
since $\odot^{1,1}$ and $\liet^{0,2}$ are in the center of the Schouten algebra and the kernel of the $\dbar$-operator,
elements in these spaces do not affect the Frobenius structure. Therefore, up to equivalence in the Gerstenhaber algebra, 
\[
\Gamma_{1}
\equiv_{\Ger}\sum_j\lambda_{j}\left(  W\wedge T_{j}\right)
+\gamma\left(W\wedge\overline{\rho}\right)
+\sum_k\alpha_{k}\left(  \overline{\rho}\wedge\oom^k\right)
\]
and hence
\begin{equation}
\Gamma
\equiv_{\Ger}\sum_j\lambda_{j}\left(  W\wedge T_{j}\right)
+\gamma\left(W\wedge\overline{\rho}\right)
+\sum_k\alpha_{k}\left(  \overline{\rho}\wedge\oom^k\right)
-\frac{1}{1-\gamma}\sum_{j,k=1}^n\lambda_j\alpha_k\psi_{jk}
\end{equation}
Let us adopt the following notations.
\begin{equation}
T=\sum_j\lambda_jT_j,
\quad
\oom=\sum_j\alpha_j\oom^j,
\quad
X=\sum_j\alpha_jT_j,
\quad
\Oomega=\sum_j\lambda_j\oom^j.
\end{equation}
In such case,
\begin{eqnarray*}
\Gamma_2 &=&-\frac12\lambda_j\alpha_k(T_j\wedge\oom^k-T_k\wedge\oom^j)=-\frac12(\lambda_jT_j\wedge \alpha_k\oom^k-\alpha_kT_k\wedge \lambda_j\oom^j)\\
&=&-\frac12(T\wedge\oom-X\wedge\Oomega).
\end{eqnarray*}
Therefore,
\begin{eqnarray}
\Gamma &\equiv_{\Ger}& W\wedge T+\gamma W\wedge\orho+\orho\wedge\oom+\frac{1}{1-\gamma}\Gamma_2
\nonumber\\
&=& W\wedge T+\gamma W\wedge\orho+\orho\wedge\oom-\frac{1}{2(1-\gamma)}(T\wedge\oom-X\wedge\Oomega).
\end{eqnarray}

We are going to compare the generators for $DGA(\dbar)$ and $DGA(\dbar_\Gamma)$, 
where
\[
DGA(\dbar)=\left( \oplus_{p,q}\lieg^{p,q}, \quad \lbra{-}{-},  \wedge, \dbar\right), 
\qquad
DGA(\dbar_\Gamma)=\left( \oplus_{p,q}\lieg^{p,q}, \quad \lbra{-}{-},  \wedge, \dbar_\Gamma\right).
\]
Recall that the structure equations for $DGA(\dbar)$ are generated by
\begin{equation}\label{dbar generator}
\dbar T_k=-\frac{i}2W\wedge\oom^k, \quad
\lbra{T_k}{\orho}=-\frac{i}2\oomega^k.
\end{equation}
Otherwise,  $\liec^{1,0}$ and $\liet^{0,1}$ are in the center of the Schouten algebra and are in the kernel
of the exterior $\dbar$-operator.

In order to compute $\dbar_\Gamma=\dbar+\ad{\Gamma}$, we take advantage of linearity:
\[
\ad{\Gamma}=\ad{W\wedge T}+\gamma\ad{W\wedge\orho}+\ad{\orho\wedge\oom}+\frac{1}{1-\gamma}\ad{\Gamma_2}.
\]
We test each summand above on $T_k$ and $\orho$ respectively.
\begin{eqnarray}
\lbra{W\wedge T}{T_k}&=&0, 
\nonumber\\
\lbra{W\wedge T}{\orho} &=& \lambda_j W\wedge\lbra{T_j}{\orho}=-\frac{i}2\lambda_jW\wedge \oom^j
=\dbar T. \label{w t orho} \\
\lbra{W\wedge\orho}{T_k}&=&-W\wedge\lbra{T_k}{\orho}=\frac{i}2W\wedge\oom^k=-\dbar T_k,
\label{w orho tk}\\
\lbra{W\wedge\orho}{\orho}&=&0.
\label{w orho orho}\\
\lbra{\orho\wedge\oom}{T_k}&=&-\oom\wedge\lbra{\orho}{T_k}=\oom\wedge\lbra{T_k}{\orho}
=-\frac{i}2\oom\wedge\oom^k.
\label{orho oom tk}\\
\lbra{\orho\wedge\oom}{\orho}&=&0.
\nonumber\\
\lbra{\Gamma_2}{T_k}&=&0. \nonumber
\end{eqnarray}
By (\ref{orho psijk}),
\begin{eqnarray}
\lbra{\Gamma_2}{\orho}&=& -\frac12 \lambda_j\alpha_k
\left(\lbra{T_j\wedge\oom^k-T_k\wedge\oom^j}{\orho}\right)
=\frac12\lambda_j\alpha_k
\left( \oom^k\wedge\lbra{T_j}{\orho}-\oom^j\wedge\lbra{T_k}{\orho}\right)\nonumber\\
&=&-\frac{i}4\lambda_j\alpha_k\left(\oom^k\wedge\oom^j-\oom^j\wedge\oom^k\right)
=-\frac{i}2\oom\wedge\Oomega. \label{gamma2 orho}
\end{eqnarray}

Now we are ready to compute $\dbar_\Gamma T_k$ and $\dbar_\Gamma\orho$. By 
(\ref{w orho tk}) and (\ref{orho oom tk}), 
\begin{eqnarray}
\dbar_\Gamma T_k&=&
\dbar T_k+\gamma \lbra{W\wedge \orho}{T_k}+\lbra{\orho\wedge\oom}{T_k}
\nonumber\\
&=&
(1-\gamma)\dbar T_k-\frac{i}2\oom\wedge\oom^k. \label{dbargamma tk}
\end{eqnarray}
On the other hand by (\ref{w t orho}) (\ref{w orho orho}) and (\ref{gamma2 orho}), 
\begin{eqnarray*}
\dbar_\Gamma\orho &=&\lbra{W\wedge T}{\orho}+\gamma\lbra{W\wedge\orho}{\orho}
+\frac{1}{1-\gamma}\lbra{\Gamma_2}{\orho}\\
&=&
\dbar T-\frac{i}{2(1-\gamma)}\oom\wedge\Oomega.
\end{eqnarray*}
Let $\sum_k\mu_kT_k+\mu \orho$ be a generic element in $\liet^{1,0}\oplus\liec^{0,1}$.  Then
\begin{eqnarray*}
&&\dbar_\Gamma(\sum_k\mu_kT_k+\mu \orho)\\
&=&(1-\gamma)\mu_k\dbar T_k-\frac{i}2\oom\wedge (\mu_k\oom^k)
+\mu \dbar T-\frac{i}2\frac{\mu}{1-\gamma}\oom\wedge\Oomega\\
&=&\dbar\left( (1-\gamma)\mu_kT_k+\mu\lambda_kT_k\right)
-\frac{i}2\oom\wedge\left( \mu_k\oom^k+\frac{\mu}{1-\gamma}\Oomega\right)\\
&=&\left( (1-\gamma)\mu_k+\mu\lambda_k\right)\dbar T_k
-\frac{i}2\frac{1}{1-\gamma}\oom\wedge
\left( (1-\gamma)\mu_k+\mu\lambda_k \right)\oom^k.
\end{eqnarray*}
It could be equal to zero only if $(1-\gamma)\mu_k+\mu\lambda_k=0$ for each $k$. i.e.
\[
\mu_k=-\frac{\mu\lambda_k}{1-\gamma}.
\]
Therefore, the only $\dbar_\Gamma$-closed element in $\liet^{1,0}\oplus\liec^{0,1}$ is
\begin{equation}
\sum_k\mu_kT_k+\mu \orho
=\mu( \orho -\sum_k\frac{\lambda_k}{1-\gamma}T_k)=
\mu (\orho-\frac{1}{1-\gamma}T).
\end{equation}
Consider a basis for $\lieg^1=\lieg^{1,0}\oplus\lieg^{0,1}$ as below.
\begin{equation}
\frac{1}{1-\gamma}T_k, \quad (1-\gamma)W+\oom, \quad
\frac{1}{1-\gamma}\oom^k, \quad \orho-\frac{1}{1-\gamma} T.
\end{equation}
The above computation shows that
\[
\dbar_\Gamma \left( (1-\gamma)W+\oom\right)=0,
\quad \dbar_\Gamma \left(\frac{1}{1-\gamma}\oom^k\right)=0,
\quad
\dbar_\Gamma\left(\orho-\frac{1}{1-\gamma} T\right)=0.
\]
Moreover, regarding the non-trivial terms by (\ref{dbargamma tk}),
\begin{eqnarray}
&&\dbar_\Gamma
\left(\frac{1}{1-\gamma}T_k\right)
\nonumber\\
 &=&\dbar T_k-\frac{i}2\oom\wedge \left(
\frac{1}{1-\gamma}\oom^k\right)
= -\frac{i}2W\wedge\oom^k-\frac{i}2\oom\wedge \left(
\frac{1}{1-\gamma}\oom^k\right)\nonumber\\
&=&-\frac{i}2\left( (1-\gamma)W+\oom\right)\wedge\left(\frac{1}{1-\gamma}\oom^k\right);
\label{dbar gamma generator}\\
&&\lbra{\frac{1}{1-\gamma}T_k}{\orho-\frac{1}{1-\gamma} T}
=\frac{1}{1-\gamma}
\lbra{T_k}{\orho}\nonumber\\
&=&\frac{1}{1-\gamma}\left(-\frac{i}2\oom^k\right)
=-\frac{i}2\left(\frac{1}{1-\gamma}\oom^k\right). \label{gamma bracket generator}
\end{eqnarray}

\begin{theorem} For any generalized deformation $\Gamma$ of the complex structure of
the primary Kodaira manifold given in Theorem \ref{Main},
there exists a natural isomorphism of Gerstenhaber algebras: 
 \[
 \left( H^\bullet_\Gamma(M), \quad \lbra{-}{-}, \quad \wedge \right) 
\quad  \cong \quad
\left( H^\bullet(M), \quad \lbra{-}{-}, \quad \wedge \right).
\]
\end{theorem}
\bproof Choose a map $\Phi$ so that $\Phi(\dbar)=\dbar_\Gamma$.
\begin{eqnarray}
&\Phi(T_k)=\frac{1}{1-\gamma}T_k, \quad \Phi(W)=(1-\gamma)W+\oom,&\\
& \Phi(\oom^k)= \frac{1}{1-\gamma}\oom^k, \quad \Phi(\orho)
=\orho-\frac{1}{1-\gamma} T.
\end{eqnarray}
Comparing the structure equations for $H^\bullet(M)$ as given in (\ref{dbar generator})
with those for $H^\bullet_\Gamma(M)$ as given in (\ref{dbar gamma generator}) and (\ref{gamma bracket generator}), 
we  complete the proof of our theorem and the main result of this paper. 
\eproof

\section{Further Development}

As noted in the Introduction, the result in Theorem \ref{main} is inspired by an investigation on the Frobenius structures on primary Kodaira surfaces. 
Although the observation in Theorem \ref{main} focuses on a generalization to primary Kodaira manifolds in all dimensions, we anticipate more because several 
key elements in this paper work in a general context. 

When working with primary Kodaira manifold, we consider it as a nilmanifold with an abelian complex structure. In this context, we make use of the fact 
\cite{MPPS-2-step} \cite[Proposition 6]{Poon1} that 
the inclusion map is a quasi-isomorphism of Gerstanhaber algebras: 
\[
\left( 
\oplus_{p, q}\lieg^{p,q},  \quad \lbra{-}{-}, \quad \wedge, \quad \dbar
\right)
\hookrightarrow
\left( 
C^\infty (M, \wedge^\bullet (\Theta\oplus \Oomega) ),  \quad \lbra{-}{-}, \quad \wedge, \quad \dbar 
\right).
\]
Given the work of others  \cite{Console} \cite{Console-Fino} \cite{MPPS-2-step} \cite{Nomizu} \cite{Rolle}, the above quasi-isomorphism is established for all nilmanifolds
with abelian complex structures \cite[Theorem 1]{CFP}. Therefore, computation of Frobenius structures could be reduced to a finite-dimensional setting
on all these manifolds. 

To make a computation effective, one may consider an \it ascending basis \rm for the algebra $\lieg^{1,0}$ and the conjugation of 
its dual  for $\lieg^{0,1}$ as in \cite{PS2}. Precisely, 
given the work of \cite{CFG} \cite{CFGU}, Salamon finds \cite{Salamon} that there exists a basis $\{ \oomega^1, \dots, \oomega^n, \oomega^{n+1} \}$ for 
$\lieg^{0,1}$ such that for all $j$, 
\begin{equation}\label{ascending}
d\omega^{j+1} \in I(\omega^1, \dots, \omega^j) \wedge I(\oomega^1, \dots, \oom^j),
\end{equation}
where $I(\omega^1, \dots, \omega^j)$ denotes the ideal generated by $\{\omega^1, \dots, \omega^j\}$. Let $\{T_1, \dots, T_{n+1}\}$ be the dual basis
for $\lieg^{1,0}$. We call both $\{T_1, \dots, T_{n+1}\}$ and $\{\omega^1, \dots, \omega^{n+1}\}$ ascending basis for the complex structure. Let $\liet^{1,0}$ and 
$\liet^{0,1}$ be the complex linear spans of $\{T_1, \dots, T_{n}\}$ and $\{\oomega^1, \dots, \oomega^{n}\}$ respectively. 
It is not hard to find that the constraints in (\ref{ascending}) imply that $\Lambda=T_{n}\wedge T_{n+1}$ is an invariant holomorphic Poisson structure \cite{PS1}. 
In \cite{PS2}, the author and his collaborator further identify conditions to secure an isomorphism between Gerstenbaher algebras:
 \[
 \left( H^\bullet_{\Lambda} (M), \quad \lbra{-}{-}, \quad \wedge \right) 
\quad  \cong \quad
\left( H^\bullet(M), \quad \lbra{-}{-}, \quad \wedge \right).
\]
On 2-step nilmanifolds with abelian complex structures, including primary Kodaira manifolds, the necessary condition could be easily stated. 
The constraint in 
(\ref{ascending}) above implies that the contraction of $d\omega^{n+1}$ is a linear map: 
\begin{equation}
d\omega^{n+1}: \liet^{1,0} \to \liet^{0,1}. 
\end{equation}
\begin{theorem}{\rm\cite[Theorem 6]{PS2}} Suppose that $M$ is a  2-step complex $(n+1)$-dimensional nilmanifold $M$ with abelian complex structure. 
Given an ascending basis with $d\omega^{n+1}$ being non-degenerate,  $\Lambda=V_{n}\wedge V_{n+1}$ is a holomorphic Poisson structure. Moreover, 
there exists a natural isomorphism between Gerstenhaber algebras:
 \[
 \left( H^\bullet_{\Lambda} (M), \quad \lbra{-}{-}, \quad \wedge \right) 
\quad  \cong \quad
\left( H^\bullet(M), \quad \lbra{-}{-}, \quad \wedge \right).
\]
\end{theorem}

\

Since the holomorphic Poisson structure represents a portion of the infinitesimal generalized deformation: 
\[
H^2(M)=H^0(M, \Theta^2) \oplus H^1(M, \Theta)\oplus H^2(M, \cal{O}),
\]
 the last theorem and the work in this paper lead to the following conjecture. 
\begin{conjecture} On a 2-step nilmanifold $M$ with abelian complex structure. Suppose that there exists an ascending basis $\{\omega^1, \dots, \omega^{n+1}\}$
such that $d\omega^{n+1}$ is non-degenerate. 
When an element $\Gamma_1 \in H^2(M)$ is integrable to  a generalized deformation $\Gamma$, 
the Gerstenbaher algebra $ \left( H^\bullet_{\Gamma} (M),  \lbra{-}{-},  \wedge \right) $ is naturally isomorphic 
to $\left( H^\bullet(M),  \lbra{-}{-},  \wedge \right).$
\end{conjecture}

\

We remark that there are examples to demonstrate that the non-degeneracy condition for $d\omega^{n+1}$ is necessary \cite{PS1} \cite{PS2}. 
On the other hand, there are many examples satisfying the conditions of this conjecture \cite{CFU} \cite{MPPS-2-step} \cite{PS1}. 

\

\noindent{\bf Acknowledgments.}  Y.S. Poon is grateful for hospitality of the Institute of Mathematical Sciences at the Chinese University of Hong Kong
during his visits in 2019.

\end{document}